\documentstyle[11pt]{article}
\title{On a stratification defined by real roots of polynomials}
\author{Vladimir Petrov Kostov \\ \\ \hspace{7cm}{\sl To Prof. Emil Horozov}} 
\date{}
\newtheorem{tm}{Theorem}

\newtheorem{rem}[tm]{Remark}
\newtheorem{rems}[tm]{Remarks}
\newtheorem{lm}[tm]{Lemma}

\newtheorem{prop}[tm]{Proposition}
\begin{document}
\maketitle

\begin{abstract}
We consider the family of polynomials $P(x,a)=x^n+a_1x^{n-1}+\ldots +a_n$, 
$x,a_i\in {\bf R}$, and the stratification of 
${\bf R}^n\cong \{ (a_1,\ldots ,a_n)|a_i\in {\bf R}\}$ defined by the 
multiplicity vector of the real roots of $P$. We prove smoothness 
of the strata and a transversality property of their 
tangent spaces.  


\end{abstract}

\section{Formulation of the result}

For $n\in {\bf N}^*$ fixed 
consider the family of polynomials $P(x,a)=x^n+a_1x^{n-1}+\ldots +a_n$, 
$x,a_i\in {\bf R}$. A {\em multiplicity vector (MV)} is a vector 
whose components are the 
multiplicities of the real roots of $P$ (for $a$ fixed) 
listed in increasing order. 
E.g. for $n=9$ the MV $[3,1,2,1]$ means that for the real 
roots $x_i$ one has $x_1=x_2=x_3<x_4<x_5=x_6<x_7$ and there is a complex 
conjugate couple about whose real and imaginary part the MV 
gives no information. 

Define the {\em length} (resp. the {\em multiplicity surplus}) of a MV 
$[r_1,r_2,\ldots ,r_q]$ as 
the integer $l=r_1+\ldots +r_q$ (resp. $r=\sum _{j=1}^q(r_j-1)=l-q$). 
The number of complex conjugate couples of roots of $P$ (counted with 
the multiplicities) equals $(n-l)/2$. 

Stratify the space 
${\bf R}^n\cong \{ (a_1,\ldots ,a_n)|a_i\in {\bf R}\}$ -- each MV 
with multiplicity surplus $r$ defines a stratum of codimension $r$. 

\begin{rems}\label{dimension}
1) It is clear that if $S_i$, $i=1,2$, are two strata  
defined by the MVs $V_i$ and if $S_1$ belongs to the closure of $S_2$, 
then for the lengths $l_i$ of the MVs $V_i$ 
one has $l_1\geq l_2$. Indeed, to pass from $S_1$ to $S_2$  
one or several real roots must bifurcate. If all new roots are real, 
then $l_1=l_2$. If there is at least one complex conjugate couple, 
then $l_1>l_2$.

2) If the MVs $V_i$ are of the same length, then $S_1$ 
belongs to the closure of $S_2$ if and 
only if $V_1$ is obtained from $V_2$ by replacing groups of 
consecutive multiplicities by their sums (call this an {\em operation of type 
A}; it corresponds to confluence 
of real roots). 

3) An {\em operation of type B} consists either in 
adding a component equal to $2$ at an 
arbitrary place of the MV (a complex conjugate couple becomes a double 
real root different from the other real roots of $P$) or in increasing of 
one of the components of the MV by $2$ 
(a complex conjugate couple becomes a double real 
root which coincides with one of the other real roots). 
If for the lengths $l_i$ of $S_i$ one has $l_1=l_2+2k$, $k\in {\bf N}$, 
then $S_1$ belongs to the closure of $S_2$ if and only if 
$V_1$ is obtained from $V_2$ either by an operation of type A  
followed by $k$ operations of type B or just by $k$ operations of type B. 

4) Given a stratum $T$ of dimension $d<n$, with MV $\vec{v}$, 
one can obtain all 
MVs of strata of dimension $d+1$ adjacent to $T$ either by replacing some 
component $m>1$ of $\vec{v}$ by two consecutive 
components $m'\geq 1,m''\geq 1$, $m'+m''=m$ or by deleting from $\vec{v}$ 
a component equal to $2$ (which means that a double real root becomes a 
complex conjugate couple). Indeed, to pass to a stratum of next dimension 
adjacent to the given one a root must bifurcate into two roots. If these 
two roots are complex conjugate, then they cannot be multiple because the 
stratification does not take into account the multiplicities of the 
complex roots.   
\end{rems}
 
The first aim of the present paper is to prove the following

\begin{tm}\label{trans}
A stratum of codimension $r$ is a smooth real contractible algebraic variety 
of dimension $n-r$. It is the graph of a smooth 
$r$-dimensional vector-function defined on the projection of the stratum 
in $Oa_1\ldots a_{n-r}$. The field of tangent spaces to the stratum is 
continuously extended to the strata belonging to its closure. The extension 
is everywhere transversal to the space $Oa_{n-r+1}\ldots a_n$. 
\end{tm}

The theorem is proved in Section~\ref{prtrans}. It is illustrated by an 
example in Section~\ref{example}. We discuss in Section~\ref{discuss} 
(in view of the theorem and of 
previous results from \cite{Ko1}) the mutual disposition of adjacent strata.

\begin{rems}
1) The theorem generalizes Theorem 1.8 from \cite{Ko2}. The latter treats the 
case when $P$ is {\em hyperbolic}, i.e. all roots are real.

2) The above stratification (or a similar one) 
has been considered (at least in some aspects) by other authors as well, 
see for instance \cite{AVG-Z}, \cite{AVaGL} and their bibliographies, 
\cite{Ka}, \cite{ShKh} and \cite{ShW}; the list is anything but exhaustive.
\end{rems}

\section{An example\protect\label{example}}

On Fig. 1. we show for $n=4$, $a_1=0$, 
the well-known picture of the swallowtail, i.e. the 
surface $\Sigma =\{ (a_2,a_3,a_4)\in {\bf R}^3|{\rm Res}(P,P')=0\}$. 
The three strata of dimension $3$ and their respective MVs are the open 
subset of ${\bf R}^3$ ``above'' $\Sigma$, with empty MV (no real roots), 
its open subset ``below'' $\Sigma$, with MV $[1,1]$, 
and the interior of the curvilinear 
pyramid $\Pi =OABC$, with MV $[1,1,1,1]$.   

The two-dimensional stratum $S$ 
defined by the MV $[2]$ is the swallowtail without the boundary of $\Pi$. 
The limits of the tangent spaces to $S$ are 
different on the different sides of the self-intersection curve $AO$ 
(excepting $O$ where the limit is unique), see Remark~\ref{remIM}. 
The strata $ABO$ and $ACO$ 
(without the boundaries) are 
defined respectively by the MVs $[2,1,1]$ and $[1,1,2]$; the stratum $BCO$ 
(without the boundary) is defined by the MV $[1,2,1]$. 

The one-dimensional strata $BO$ and $CO$ (without the point $O$) 
are defined respectively by the MVs $[3,1]$ and 
$[1,3]$ while $AO$ (without the point $O$) is defined by the MV $[2,2]$. 
Finally, the point $O$ is defined by the MV $[4]$. 

The self-intersection curve $AO$ has an analytic continuation $OD$ 
(represented by a dashed line on Fig. 1) which consists of polynomials of 
the form $(x^2+b^2)^2$, $b\in {\bf R}$, i.e. having a conjugate complex 
couple of roots of multiplicity $2$. This curve is not a stratum of the 
stratification (although it belongs to $\Sigma$) because the multiplicities 
of the complex roots are not taken into account, see part 4) of 
Remarks~\ref{dimension}. 

\begin{rem}\label{remIM}
The following proposition is proved in \cite{Me}, p. 52-53.

\begin{prop}\label{propIM}
If $P$, $P_1$, $\ldots$, $P_r$ are   
monic polynomials where $P=P_1P_2\ldots P_r$ and $P_i$ have two by two 
no root in common, 
then there exist neighbourhoods $U$, $U_1$, $\ldots$, $U_r$ of 
$P$, $P_1$, $\ldots$, $P_r$ such 
that the product map $U_1\times \ldots \times U_2\rightarrow U$, 
$(Q_1,\ldots ,Q_r)\mapsto Q$, is 
a diffeomorphism. 
\end{prop}

This result implies that at any point of the open arc $AO$ the swallowtail 
is locally diffeomorphic to two intersecting planes in ${\bf R}^3$. Indeed, 
on $AO$ the polynomial $P$ has two double real roots, hence, is of the form 
$(x-a)^2(x+a)^2$, $a\in {\bf R}$. One can set $P_1=(x-a)^2$, $P_2=(x+a)^2$ 
and $Q_1=(x-a)^2+\alpha (x-a)+\beta$, $Q_2=(x+a)^2-\alpha (x+a)+\gamma$. The 
two planes are defined respectively by the variables $(\alpha ,\beta )$ and 
$(\alpha ,\gamma )$. 

Proposition~\ref{propIM} can be applied in the general case as well 
(i.e. when $P$ has any MV) to understand what the set 
$\{ (a_1,\ldots ,a_n)\in {\bf R}^n|{\rm Res}(P,P')=0\}$ is locally 
like up to a diffeomorphism.
\end{rem}

\section{On the mutual disposition of strata\protect\label{discuss}}

Consider a point $A$ of a stratum $U$ of dimension $s\leq n-2$. Intersect 
$U$ by the affine space ${\cal F}$ 
of dimension $2$ containing $A$ and parallel to 
$Oa_{s+1}a_{s+2}$. By Theorem~\ref{trans} the intersection is the point $A$. 
The intersections with ${\cal F}$ of the strata of dimension $s+1$  
are curves containing $A$ in their closures and having 
non-vertical limits at $A$ of their tangent lines (Theorem~\ref{trans}). 
Each such curve (considered locally, at $A$) projects only ``to the left'' 
or ``to the right'' of $A$ on $Oa_{s+1}$. The intersections 
with ${\cal F}$ of the strata of dimension $s+2$  
are sectors delimited by these curves.

It is explained in \cite{Ko1}, Subsection 1.2, how the above curves are 
situated near $A$ in the case when $U$ is a stratum of hyperbolic polynomials. 
We generalize here these results for the case of arbitrary stratum $U$. 

Denote by $\vec{v}=[r_1,\ldots ,r_q]$ the MV of the stratum $U$. Denote by 
$U_{i,j}$ the stratum with MV obtained from $\vec{v}$ by replacing the 
component $r_i$ by two components -- $j$, $r_i-j$ -- where 
$j=1,\ldots ,r_i-1$. If $r_i=2$, denote by $V_i$ the stratum whose MV is 
obtained from $\vec{v}$ by deleting the component $r_i$. By part 4) of 
Remarks~\ref{dimension}, these are all strata of dimension $s+1$ adjacent to 
$U$. We use the notation $U_{i,j}$, $V_i$ also for the intersections of 
these strata with ${\cal F}$. In what follows we assume that the equations 
of the limits at $A$ of the tangent lines to the curves $U_{i,j}$ are given 
in the form $a_{s+2}=k_ia_{s+1}+\theta _i$.
 
\begin{lm}\label{slope}
The slopes of the limits at $A$ of the tangent lines to the curves 
$U_{i_1,j_1}$, $U_{i_2,j_2}$, are the same for $i_1=i_2$ and different for 
$i_1\neq i_2$.
\end{lm}

The lemma is proved by full analogy with Lemma 16 from \cite{Ko1} (using 
Proposition~\ref{propIM}). In what follows we assume that the equations 
of the limits at $A$ of the tangent lines to the curves $U_{i,j}$ are given 
in the form $a_{s+2}=k_ia_{s+1}+\theta _i$. Here $k_i$ is the slope of the 
limit of the tangent line.

\begin{lm}\label{UV}
If $r_i=2$, then the tangent lines to $U_{i,1}$ and $V_i$ are the same. 
The curves $U_{i,1}$, $V_i$ and the point $A$ are parts of one and the same 
curve smooth at $A$.
\end{lm}

{\em Proof:}

In the particular case when $P=x^2+\lambda$ the strata $U_{i,1}$ and $V_i$ are 
the half-lines $\{ \lambda <0\}$ and $\{ \lambda >0\}$. In the general case 
the lemma 
is proved by analogy with Lemma 16 from \cite{Ko1} (using 
Proposition~\ref{propIM}).~~~~$\Box$   
 
\begin{lm}\label{slopebis}
For the slopes $k_i$ of the limits at $A$ of the tangent lines to 
the curves $U_{i,j}$ one has $k_1>\ldots >k_q$. 
\end{lm}

{\em Proof:}

$1^0$. In the case when the stratum $U$ consists of hyperbolic polynomials 
the lemma 
is proved in \cite{Ko1}, see Lemma~22 there. Suppose that $U$ does not 
consist of hyperbolic polynomials. Denote by $U'$ the stratum whose MV is 
obtained from $\vec{v}$ by adding to the right $n-l$ components equal to $1$. 
Hence, $U'$ consists of hyperbolic polynomials. 

$2^0$. Choose a point $B\in U'$. Connect it with $A$ by a continuous curve 
(parametrized by $\sigma \in [0,1]$) passing 
only through strata with MVs whose first $q$ components are the same as the 
ones of $\vec{v}$ and such that for $(n-l)/2$ distinct values $\sigma _{\nu }$ 
of $\sigma$ the 
greatest two of the real roots become equal after what they form a complex 
conjugate couple. The slopes $k_i$ can be defined for any $\sigma \in [0,1]$. 
For different $i$ they remain different throughout the deformation 
(even for $\sigma =\sigma _{\nu }$ which can be proved like Lemma~\ref{slope}).
For all $\sigma$ they are finite. Hence, their order is the same at 
$A$ and at $B$.~~~~$\Box$

\begin{lm}\label{leftright}
If $i$ is even, then the projection on $Oa_{s+1}$ of $U_{i,j}$ is ``to the 
right'' of the one of $A$; if $i$ is odd, then ``it is on its left''. 
\end{lm}

{\em Proof:}

In the case when $U$ consists of hyperbolic polynomials this is Lemma~18 
from \cite{Ko1}. In the general case the lemma is proved like the previous one 
-- being ``to the left or to the right'' does not change throughout the 
deformation, because it depends continuously on $\sigma$, 
i.e. in fact it does not depend on $\sigma$.~~~~$\Box$

\begin{rem}
It follows from Lemmas~\ref{UV} and \ref{leftright} that if $r_i=2$ and 
if the projection of $U_{i,1}$ on $Oa_{s+1}$ is ``to the right'' (resp. ``to 
the left'') of the one of $A$, then the projection of $V_i$ 
on $Oa_{s+1}$ is ``to the left'' (resp. ``to 
the right'') of the one of $A$.
\end{rem}

\begin{lm}\label{updown}
For $i$ fixed the curve $U_{i,j_1}$ is ``above'' the curve $U_{i,j_2}$ 
if and only if either $i$ is odd and $j_1>j_2$ or $i$ is even and $j_1<j_2$.
\end{lm}

{\em Proof:}

If the stratum $U$ consists of hyperbolic polynomials, then this is Lemma~20 
from \cite{Ko1}. If not, then use the same deformation as in the proof of 
Lemma~\ref{slopebis}. For any value of $\sigma$ and for any $i$ fixed 
the curves $U_{i,j}$ have the same mutual disposition. Indeed, one can apply 
Proposition~\ref{propIM} -- the curves $U_{i,j}$ with one and the same $i$ 
correspond to one and the same neighbourhood $U_i$ and to the respective 
curves constructed after the polynomial $(x-a)^{r_i}$. The latter's 
mutual disposition does not depend on $\sigma$ or on $a$.~~~~$\Box$

\begin{rem}
On Fig. 2 we show the curves $U_{i,j}$ and $V_i$ in the case when $l=10$, 
$n=10+2h$, $h\in {\bf N}$. The reader can check the above lemmas on this 
example.
\end{rem}
   
\begin{rem} 
If for some $i$ one has $r_i=r_{i+1}=2$, then the MVs of the strata 
$V_i$ and $V_{i+1}$ are the same. Hence, this is one and the same stratum but 
it gives rise to two (or more) different curves $V_i$ at $A$, with different 
slopes of their tangent lines at $A$ Lemma~\ref{slope}). The simplest 
example of such a situation 
is the one of the previous section (the stratum defined by the 
MV $[2]$ admitting 
two different limits of the tangent space along the curve $AO$). 
\end{rem}

\section{Proof of Theorem~\protect\ref{trans}\protect\label{prtrans}}

$1^0$.  
Contractibility follows from the 
contractibility of the parameter 
space (the roots play the role of parameters and 
define the coefficients $a_i$ via the Vieta formulas; these formulas define 
a homeomorphism). 

Further we prove some statements of the 
theorem not for the stratification of ${\bf R}^n$, the space of the 
coefficients $a_i=(-1)^i\sigma _i$ (where $\sigma _i$ 
is the $i$-th symmetric function of the roots of $P$ counted with the 
multiplicities), but for the space of the Newton 
functions $b_i$ which are the sums of the $i$-th powers of the roots. The 
statements formulated for the two spaces (of the quantities $a_i$ or $b_i$) 
are equivalent because 
there exist polynomials $q_j$, $q^*_j$ such that 

\[ ja_j=-nb_j+q_j(b_1,\ldots ,b_{j-1})~~,~~
nb_j=-ja_j+q^*_j(a_1,\ldots ,a_{j-1})~.\]

$2^0$. Denote by $(x_1,\ldots ,x_{n-r})$ the roots of $P$ where the real 
roots are distinct while the complex ones might not be all distinct. 
The Jacobian matrix corresponding to the 
mapping $(x_1,\ldots ,x_{n-r})\mapsto (b_1,\ldots ,b_{n-r})$ 
is obtained from the Vandermonde matrix $W(x_1,\ldots ,x_{n-r})$ by 
multiplying the columns of the real roots by their respective multiplicities. 
Hence, if all complex roots are distinct, then the determinant $J$ 
of the matrix $W$ is nonzero and at such a point the mapping is a local 
diffeomorphism which means that the stratum is locally of dimension $n-r$. At 
such a point one can express the roots $x_i$ as smooth functions of 
$(b_1,\ldots ,b_{n-r})$, and then express $(b_{n-r+1},\ldots ,b_n)$ as 
smooth functions of $x_i$; hence, as smooth functions of 
$(b_1,\ldots ,b_{n-r})$. 

$3^0$. Prove the smoothness of the stratum regardless of whether the complex 
roots are all distinct or not. (We use the same ideas here as in \cite{Me}, 
p. 52-53.) To this end set $P=QR$, 
$Q=x^{2k}+c_1x^{2k-1}+\ldots +c_{2k}$, 
$R=x^{n-2k}+d_1x^{n-2k-1}+\ldots +d_{n-2k}$ where all roots of $Q$ are 
complex and all roots of $R$ are real. The mapping 

\[ (c_1,\ldots ,c_{2k},d_1,\ldots ,d_{n-2k})\mapsto (a_1,\ldots ,a_n)\]
(where $c_i$, $d_j$ are regarded as free parameters) is a local 
diffeomorphism. Indeed, its Jacobian matrix is the Sylvester matrix of $Q$ 
and $R$; the latter's determinant equals Res$(Q,R)$ which is nonzero because 
$Q$ and $R$ have no root in common. 

$4^0$. {\em The field of tangent spaces to the given stratum is continuously 
extended to the strata of lower dimension belonging to its closure and to the 
points where some complex roots coincide. The extension is everywhere 
transversal to the space $Oa_{n-r+1}\ldots a_n$.}

The rest of the theorem follows from the statement. The latter implies in 
particular that even in a neighbourhood of a point of the stratum 
where some complex roots coincide, the stratum is locally the graph of a 
smooth $r$-dimensional function defined on the projection of the stratum in 
$Oa_1\ldots a_{n-r}$, see $2^0-3^0$. 

$5^0$. To prove the statement from $4^0$ compute the partial derivatives 
$\partial b_k/\partial b_u$, $k\geq n-r+1$, $u\leq n-r$ bearing in mind that 
$b_j$ is the $j$-th Newton function of the roots $x_i$. (We follow here 
the same ideas as the ones used in the proof of Theorem 1.8 from \cite{Ko2}.) 
Denote by $m_i$ the 
quantity equal to the multiplicity of $x_i$ if $x_i$ is a real root and to $1$ 
if it is a complex one. One has 

\begin{equation}\label{partial} 
\partial b_k/\partial b_u=
\sum _{i=1}^{n-r}(\partial b_k/\partial x_i)(\partial x_i/\partial b_u)=
k\sum _{i=1}^{n-r}(m_ix_i^{k-1})(\partial x_i/\partial b_u)=
k\sum _{i=1}^{n-r}(m_ix_i^{k-1}A_{u,i})/w
\end{equation}
where $w=$det$\| \partial b_j/\partial x_{\nu }\| =g\prod _{q<v}(x_q-x_v)$, 
$g\neq 0$, and $A_{u,i}$ is the cofactor of the element 
$\partial b_u/\partial x_i$ in the matrix 
$\| \partial b_j/\partial x_{\nu }\|$. 

$6^0$. Put $x_{\mu}=x_{\nu}$ in (\ref{partial}). 
Then for $i\neq \mu ,\nu$ one has 
$A_{u,i}=0$ (two proportional columns). Suppose first that the roots 
$x_{\mu}$, $x_{\nu}$ are complex. Then one has $m_{\mu }=m_{\nu }=1$ and 
$m_{\mu}A_{u,\mu}+m_{\nu}A_{u,\nu}=A_{u,\mu}+A_{u,\nu}=0$ (to be checked 
directly). This means that the 
numerator of the right hand-side of (\ref{partial}) is $0$ when 
$x_{\mu}=x_{\nu}$, i.e. it is representable in the form 
$wh(x_1,\ldots ,x_{n-r})$ for some polynomial $h$. Hence, 

\begin{equation}\label{partial1}
\partial b_k/\partial b_u=h
\end{equation}   

If the roots $x_{\mu}$, $x_{\nu}$ are real, then one again 
checks directly that 
$m_{\mu}A_{u,\mu}+m_{\nu}A_{u,\nu}=0$ and again one has (\ref{partial1}). 

$7^0$. The closure of the stratum can be defined by a continuous 
parametrization of 
the roots $x_i$ by some parameters $z$ (one can choose as such parameters 
part of the variables $b_j$; in general, these variables are more than the 
parameters needed). By (\ref{partial1}), the partial 
derivatives $\partial b_k/\partial b_u$ are bounded continuous 
functions of the parameters $z$. Hence, the limits of these partial 
derivatives exist on the closure of the stratum. This proves the statement.

The theorem is proved. ~~~~$\Box$

\end{document}